\documentclass[12pt]{amsart}
\usepackage{amsfonts,amsmath,amssymb,amscd,amsthm,enumerate,tikz}
\usepackage{amsrefs}
\usepackage{geometry}
\usepackage{adjustbox}
\usepackage[all]{xy}
\usepackage{stmaryrd}
\usetikzlibrary{arrows.meta,decorations.markings}
\newtheorem{theorem}{Theorem}
\newtheorem{prop}[theorem]{Proposition}

\newtheorem{cor}[theorem]{Corollary}
\theoremstyle{definition}

\newtheorem{rem}[theorem]{Remark}
\newtheorem{mydef}[theorem]{Definition}
\newtheorem{example}[theorem]{Example}

\renewcommand{\epsilon}{\varepsilon}

\def\<{\langle}
\def\>{\rangle}

\usepackage{multicol} %%%%%%%%%%%

\usepackage{tikz}
\usepackage{caption}
\usepackage{subcaption}

\begin{document}
\title[Lower bounds for hypergraph independence]{An Efficiently Computable Lower Bound for the Independence Number of Hypergraphs}
\author[M.\ Aldi, T.\ Gabrielsen, D.\ Grandini, J.\ Harris, K.\ Kelley]{Marco Aldi, Thor Gabrielsen, Daniele Grandini, Joy Harris, Kyle Kelley}

\begin{abstract}
Let $k\ge2$ be fixed. We study the integer lower bound
$\ell(G)$ obtained by inverting the classical counting
inequality for Tur\'an systems. For a $k$-uniform hypergraph
with $n$ vertices and $m$ edges, the bound can be evaluated
exactly by binary search in time polynomial in the binary
lengths of $n$ and $m$. We exhibit separations from the Tur\'an-Spencer and
Caro-Tuza bounds for every fixed $k\ge3$, and from the
Csaba-Plick--hokoufandeh bound in the $3$-uniform case.
For the Caro-Tuza comparison, the separation grows linearly in $k$ on infinitely many regular
$k$-uniform hypergraphs.
\end{abstract}
\maketitle

\section{Introduction} 
In a hypergraph $G$, an {\it independent set} is a subset of the vertices that contains no edges. The cardinality of the largest independent set is the independence number, denoted $\alpha(G)$. Because calculating $\alpha(G)$ is NP-hard in general, and approximating $\alpha(G)$ within a factor of $n^{1-\epsilon}$ for every fixed $0<\epsilon<1$ is NP-hard for explicitly represented $n$-vertex $k$-uniform hypergraphs for each fixed $k\ge2$ \cites{Haastad1999,Zuckerman2007}, efficiently computable lower bounds for $\alpha(G)$ are essential for both theoretical and practical applications.

Existing lower bounds for $\alpha(G)$ generally fall into two categories based on their required input information. The first category consists of global parameter bounds that rely solely on summary parameters of the hypergraph, such as the number of vertices $n$, edges $m$, and uniformity $k$. The classical example is the Tur\'an–Spencer bound \cites{Spencer1972, Eustis2013}. These bounds require no structural graph information beyond $n, m,$ and $k$. The second category consists of local degree-based bounds that incorporate fine-grained structural data, specifically the degree sequence $(d(v))_{v \in V(G)}$. Prominent examples include the Caro–Tuza bound  \cite{CaroTuza1991} and the Csaba–Plick–Shokoufandeh bound  \cite{CsabaPlickShokoufandeh2012}. While degree-based bounds utilize local degree distributions, they do not strictly dominate global parameter bounds across all hypergraphs; as we demonstrate in Section 3, $\ell(G)$ can yield strictly higher lower bounds than both Caro–Tuza and Csaba–Plick–Shokoufandeh for certain families of hypergraphs.

The bound studied in this paper belongs to the classical theory of
Tur\'an systems and covering designs; see \cite{Sidorenko1995}.
If a $k$-uniform hypergraph $G$ has at least one edge and
$\alpha(G)=a$, then its edge set is a Tur\'an $(n,a+1,k)$-system:
every $(a+1)$-subset of its vertex set contains an edge. The
elementary counting bound
\[
T(n,s,k)\ge \frac{\binom{n}{k}}{\binom{s}{k}}
\]
therefore gives a lower bound for $\alpha(G)$ by inversion.
We denote this integer lower bound by $\ell(G)$, with the convention
$\ell(G)=n$ when $G$ is edgeless. 

Our contribution concerns the evaluation and comparison of this
explicit bound. For fixed uniformity $k$, we give an exact
binary-search procedure and analyze its bit complexity. We prove
that $\ell(G)$ never exceeds the rounded average-degree Tur\'an
bound for graphs. In contrast, for every $k\ge3$, we identify
parameter ranges containing infinitely many $k$-uniform hypergraphs
on which $\ell(G)$ exceeds the Tur\'an-Spencer bound, and we
construct families illustrating how the separation between $\ell(G)$ and the
Caro-Tuza lower bound can grow linearly with $k$. In the $3$-uniform case, we also give examples of separation between $\ell(G)$ and the 
Csaba-Plick-Shokoufandeh lower bound.
These comparisons concern the specified bounds, rather than
optimality among all lower bounds determined by $(n,m,k)$.
In Remark~\ref{rem:covering-comparison}, we explain how the
counting estimate underlying $\ell(G)$ relates to classical
recursive covering bounds and to de Caen's lower bound. 

Our interest in $\ell(G)$ originates from our work on independence theory of 2-step nilpotent Lie algebras \cite{AGGHK}. The cohomological interpretation of $\alpha(G)$ given in \cite{AGGHK}, together with the generalization of the Dani-Mainkar construction to hypergraphs \cite{AldiBevins}, naturally yield an alternate approach to the bounds which in the present paper we describe from a purely combinatorial perspective.

\section{The Main Theorem}

In this section, we introduce the lower bound $\ell(G)$ for the independence number of an arbitrary $k$-uniform hypergraph and establish its basic properties, including its exactness for complete  hypergraphs. 

\bigskip

Throughout this paper, a hypergraph $G = (V, E)$ is assumed to be finite and simple (containing no repeated edges). We denote $n = \vert{}V(G)\vert{}$ and $m = \vert{}E(G)\vert{}$. For a $k$-uniform hypergraph, we assume $n \ge k \ge 2$ unless explicitly stated otherwise.

\begin{mydef}
For integers $n\ge s\ge k$, a \emph{Tur\'an $(n,s,k)$-system}
is a family $\mathcal{F}\subseteq \binom{[n]}{k}$ such that every
$s$-subset of $[n]$ contains at least one member of $\mathcal{F}$.
We denote by $T(n,s,k)$ the minimum cardinality of such a family.
We use the natural convention
\[
T(n,k,k)=\binom{n}{k}.
\]
\end{mydef}

Equivalently, by taking complements,
\[
T(n,s,k)=C(n,n-k,n-s),
\]
where $C(v,b,t)$ denotes the covering number: the minimum number of
$b$-subsets of a $v$-set needed to contain every $t$-subset.

\begin{prop}\label{prop:turan-interpretation}
Let $G$ be a finite, simple $k$-uniform hypergraph with $n$ vertices
and $m\ge 1$ edges, and let $a=\alpha(G)$. Then
\[
m\ge T(n,a+1,k).
\]
In particular,
\begin{equation}\label{eq:counting-turan}
\binom{n}{a+1}
\le
m\binom{n-k}{a-k+1},
\end{equation}
or, equivalently,
\begin{equation}\label{eq:counting-turan-equivalent}
m\binom{a+1}{k}\ge \binom{n}{k}.
\end{equation}
\end{prop}

\begin{proof}
By the definition of $a=\alpha(G)$, every $(a+1)$-subset of $V(G)$
contains an edge of $G$. Hence $E(G)$ is a Tur\'an
$(n,a+1,k)$-system, proving $m\ge T(n,a+1,k)$.

Consider pairs $(e,S)$ such that
$e\in E(G)$, $S\in\binom{V(G)}{a+1}$, and $e\subseteq S$. Every
$(a+1)$-set occurs in at least one such pair, whereas every fixed
edge occurs in exactly
\[
\binom{n-k}{a-k+1}
\]
pairs. This proves \eqref{eq:counting-turan}. The equivalent form
\eqref{eq:counting-turan-equivalent} follows from
\[
\frac{\binom{n}{a+1}}
     {\binom{n-k}{a-k+1}}
=
\frac{\binom{n}{k}}
     {\binom{a+1}{k}}.
\]
\end{proof}

\begin{mydef}\label{def:counting-bound}
Let $G$ be a finite, simple $k$-uniform hypergraph with $n$ vertices
and $m$ edges. If $m=0$, set
\[
\ell(G)=n.
\]
If $m\ge 1$, set
\begin{align}\label{eq:counting-bound}
\ell(G)
&=
\min\left\{
i\in\{k-1,\ldots,n-1\}\,\mid \,
m\binom{i+1}{k}\ge\binom{n}{k}
\right\}\\
&= \min\left\{
i\in\{k-1,\ldots,n-1\}\,\mid \, f(i)\le m
\right\}
\end{align}
where
\begin{equation}\label{eq:3}
f(i) = \frac{\binom{n}{k}}{\binom{i+1}{k}}=\frac{n(n-1)\cdots(n-k+1)}{(i+1)i\cdots (i-k+2)}
\end{equation}
\end{mydef}

\begin{theorem}\label{thm:counting-bound}
For every finite, simple $k$-uniform hypergraph $G$, $\ell(G)\le \alpha(G)$.
\end{theorem}

\begin{proof}
The assertion is immediate when $m=0$. Suppose that $m\ge1$ and put
$a=\alpha(G)$. Proposition \ref{prop:turan-interpretation} gives
\[
m\binom{a+1}{k}\ge\binom{n}{k}.
\]
Therefore $a$ belongs to the set in
\eqref{eq:counting-bound}, and its minimum is at most $a$.
\end{proof}

\begin{rem}
The bound $\ell(G)$ depends only on the integer parameters $(n, m, k)$ rather than the explicit edge set of $G$. Because, for fixed $k\ge 2$, $f(i)$ as in \eqref{eq:3}
is strictly decreasing on $\{k-1, \ldots, n\}$, the minimal index $i$ satisfying $f(i) \le m$ can be located using binary search. The input consists of $n$ and $m$, encoded in binary, with $k$
fixed. If $m=0$, return $n$. Otherwise, binary search on
$\{k-1,\ldots,n-1\}$ uses exact integer comparisons of
\[
m\prod_{j=0}^{k-1}(i+1-j)
\ge
\prod_{j=0}^{k-1}(n-j).
\]
Since the binary search requires at most \(\lceil\log_2(n-k+2)\rceil\) comparisons, the total arithmetic complexity of calculating $\ell(G)$ is $O(\log n)$ operations. In a bit complexity model with $B = O(\log n)$ bit-length integers, $\ell(G)$ is computable in $O( \log n \cdot M(\log n))$ time, where $M(B)$ is the complexity of $B$-bit integer multiplication. This establishes that $\ell(G)$ is computable in polynomial time in the binary input length $O(\log n + \log(m+1))$. 
\end{rem}

\begin{rem}\label{rem:hardness}
For every fixed integer $k\ge2$ and every fixed
$0<\varepsilon<1$, approximating the maximum independent set
within a factor $n^{1-\varepsilon}$ is NP-hard for finite simple
$n$-vertex $k$-uniform hypergraphs given explicitly by their
vertex and edge lists. For graphs, this follows from
Zuckerman's inapproximability theorem for Max Clique
\cite{Zuckerman2007}*{Theorem~1.1}, which derandomizes
H\aa stad's result \cite{Haastad1999}, by taking complements,
since $\alpha(G)=\omega(\overline G)$. The extension to each
fixed $k\ge3$ follows from the following elementary reduction. Given a graph $G$ on $n\ge k$ vertices, let $H_k(G)$ be the
$k$-uniform hypergraph on the same vertex set whose edges are
all $k$-subsets containing an edge of $G$. This construction
takes polynomial time for fixed $k$. Every independent set
of $G$ is independent in $H_k(G)$, and every independent set
of $H_k(G)$ of size at least $k$ is independent in $G$.
Consequently,
\[
\alpha(H_k(G))=\max\{\alpha(G),k-1\}.
\]
To handle the small-set exception, enumerate the independent
sets of $G$ of size at most $k-1$. If no independent
$(k-1)$-set exists, this determines $\alpha(G)$ exactly.
Otherwise, fix such a set $J$. Then
$\alpha(H_k(G))=\alpha(G)$, and an independent set returned
for $H_k(G)$ can be used unchanged when its size is at least
$k$, or replaced by $J$ otherwise. Thus an approximation for
$H_k(G)$ yields an approximation for $G$ with the same
factor, without changing the number of vertices.
\end{rem}

\begin{rem}
\label{rem:covering-comparison}
For integers $n\ge s\ge k\ge2$, write
\[
B(n,s,k)
=
\frac{\binom{n}{k}}{\binom{s}{k}}=f(s-1).
\]
Thus, for a hypergraph with $m\ge1$ edges, $\ell(G)$ is the
least integer $i\in\{k-1,\ldots,n-1\}$ such that
\[
B(n,i+1,k)\le m\,.
\]
In our notation, the inequality
of Sch\"onheim and, independently, Katona, Nemetz and
Simonovits \cites{Schonheim1964,KNS} is
\begin{equation}\label{eq:covering-recursion}
T(n,s,k)
\ge
\left\lceil
\frac{n}{n-k}\,T(n-1,s,k)
\right\rceil,
\qquad n>s\ge k.
\end{equation}
Starting with $T(s,s,k)=1$, iterating
\eqref{eq:covering-recursion}, and omitting the ceiling
operations gives
\[
T(n,s,k)
\ge
\prod_{j=s+1}^{n}\frac{j}{j-k}
=
\frac{\binom{n}{k}}{\binom{s}{k}}
=
B(n,s,k).
\]
Retaining the intermediate ceilings can give a stronger estimate.
Consequently, $\ell(G)$ should be viewed simply as an inversion of
the elementary counting bound. Another classical estimate is de Caen's inequality
\cite{deCaen1983}; see also \cite{Sidorenko1995}:
\begin{equation}\label{eq:decaen-covering}
T(n,s,k)
\ge
B_{\rm dC}(n,s,k)
:=
\frac{n-s+1}{n-k+1}
\frac{\binom{n}{k}}{\binom{s-1}{k-1}}.
\end{equation}
A direct comparison gives
\[
\frac{B_{\rm dC}(n,s,k)}{B(n,s,k)}
=
\frac{s(n-s+1)}{k(n-k+1)}
=
1+
\frac{(s-k)(n+1-s-k)}{k(n-k+1)}.
\]
In particular, de Caen's estimate is strictly stronger than
the elementary count whenever
\[
k<s<n+1-k.
\]
The two unrounded estimates are not uniformly ordered over
the full parameter range. This difference can lead to a strictly stronger lower bound
for the independence number. 

For a nonempty $k$-uniform hypergraph $G$, define
\[
\ell_{\rm dC}(G)
=
\min\left\{
i\in\{k-1,\ldots,n-1\}:
B_{\rm dC}(n,i+1,k)\le m
\right\}.
\]
For an edgeless hypergraph, set $\ell_{\rm dC}(G)=n$;
then $\ell_{\rm dC}(G)=\alpha(G)$.
If $m\ge1$, Proposition~\ref{prop:turan-interpretation}
and \eqref{eq:decaen-covering} give
\[
m\ge T(n,\alpha(G)+1,k)
  \ge B_{\rm dC}(n,\alpha(G)+1,k).
\]
Hence $\ell_{\rm dC}(G)\le\alpha(G)$ in every case. For example, let $G$ be any
simple $3$-uniform hypergraph with $n=10$ vertices and $m=32$
edges. The definition gives $\ell(G)=3$, since
\[
32\binom{3}{3}
<
\binom{10}{3}
\le
32\binom{4}{3}.
\]
On the other hand, \eqref{eq:decaen-covering} yields
\[
T(10,4,3)
\ge
\frac{7}{8}\frac{\binom{10}{3}}{\binom{3}{2}}
=
35.
\]
If $\alpha(G)\le3$, every $4$-subset of $V(G)$ would contain
an edge, forcing
\[
32=|E(G)|\ge T(10,4,3)\ge35,
\]
a contradiction. Hence $\alpha(G)\ge4$, whereas the
elementary counting inversion gives only $\ell(G)=3$. On the other hand, consider the 4-uniform hypergraph with $V(G)=\{1,2,\ldots,7\}$ and $E(G)=\{\{1,2,3,4\},\{1,5,6,7\}\}$. Since $f(5)=\frac{7}{3}>2$ and $f(6)=1\le 2$, we conclude that $\ell(G)=6=\alpha(G)$. At the same time
\[
B_{\rm dC}(7,6,4)=\frac{7}{4}<2<\frac{105}{16}=B_{\rm dC}(7,5,4)
\]
implies that the algebraic inversion of the de Caen estimate
only gives the lower bound $\alpha(G)\ge \ell_{\rm dC}(G)=5$. Thus neither of the two integer bounds dominates the other.
If $r=\ell_{\rm dC}(G)<n-k+1$, the preceding ratio identity
gives
\[
f(r)\le B_{\rm dC}(n,r+1,k)\le m,
\]
and hence $\ell(G)\le\ell_{\rm dC}(G)$.

The separating families considered below also yield
separations for $\ell_{\rm dC}$. Indeed, for every nonempty,
noncomplete $k$-uniform hypergraph,
\[
B_{\rm dC}(n,k,k)=\binom nk>m,
\]
so
\begin{equation}\label{eq:dC}
\ell_{\rm dC}(G)\ge k\,.
\end{equation}

For a broader account of Tur\'an-system bounds, see
\cite{Sidorenko1995}; for recent work on asymptotic upper
bounds for the minimum size of Tur\'an systems, see
\cite{LiuPikhurko2026}.
\end{rem}

\begin{cor}\label{cor:3}
Let $G$ be a finite simple graph with $n\ge 2$ vertices and $m$ edges.
If $m=0$, then
\[
\ell(G)=\alpha(G)=n.
\]
If $m\ge 1$, then
\begin{equation}\label{eq:graph-counting}
\alpha(G)\ge \ell(G)
=
\left\lceil
\frac{\sqrt{m^2+4mn(n-1)}-m}{2m}
\right\rceil .
\end{equation}
Equivalently,
\[
\alpha(G)\ge
\left\lceil
\frac{\sqrt{1+\frac{4n(n-1)}{m}}-1}{2}
\right\rceil .
\]
\end{cor}

\begin{proof}
The assertion is immediate when $m=0$, since $G$ is then edgeless. Suppose that $m\ge 1$. Setting $k=2$ in the definition, we obtain that $\ell(G)$ is the minimum of the $i\in\{1,\ldots,n-1\}$ such that $f(i)\le m$ i.e.\
\[
mi(i+1)\ge n(n-1),
\]
or
\[
mi^2+mi-n(n-1)\ge 0.
\]
The positive root of the corresponding quadratic equation is
\[
\rho
=
\frac{\sqrt{m^2+4mn(n-1)}-m}{2m}.
\]
Since the quadratic is strictly increasing for $i\ge0$, the least
integer satisfying the inequality is $\lceil\rho\rceil$. Therefore,
\[
\ell(G)
=
\left\lceil
\frac{\sqrt{m^2+4mn(n-1)}-m}{2m}
\right\rceil\, .
\]
\end{proof}

\begin{rem}\label{rem:5}
Let $G$ be a simple $k$-uniform hypergraph with $n$ vertices and $m$ edges. The definition of $\ell(G)$ can be restated explicitly as follows:
\begin{enumerate}[(i)]
    \item For $r \in \{k, \ldots, n-1\}$, $\ell(G) = r$ if and only if
    \begin{equation}\label{eq:7}
    f(r) \le m < f(r-1)\,.
    \end{equation}
    \item For $r = k - 1$, $\ell(G) = k - 1$ if and only if $m \ge f(k-1) = \binom{n}{k}$.
\end{enumerate}
Alternatively, adopting the convention $f(k-2)=\infty$, for
$m\ge1$ the characterization \eqref{eq:7} holds uniformly
for all $r\in\{k-1,\ldots,n-1\}$. The remaining case is
$\ell(G)=n$ if and only if $m=0$.
\end{rem}

\begin{prop}\label{prop:5} Let $G$ be a $k$-uniform hypergraph on $n$ vertices. Then:  
\begin{enumerate}[1)]
\item $\ell(G) = k-1$ if and only if $G$ is complete.  
\item $\ell(G) = n$ if and only if $G$ is edgeless ($m=0$).  
\end{enumerate}
In both cases, $\ell(G) = \alpha(G)$.

\end{prop}

\begin{proof}
If $\ell(G) = k - 1$, then by definition of $\ell(G)$, we must have $f(k-1) \le m$. Evaluating $f$ at $i = k - 1$ yields:
\begin{equation*}
f(k-1) = \frac{\binom{n}{k}}{\binom{k}{k}} = \binom{n}{k} \le m\,.
\end{equation*}
Since $G$ is simple, it contains at most $\binom{n}{k}$ edges. Thus, $m = \binom{n}{k}$, which implies $G$ is complete. Conversely, if $G$ is complete, its independence number is $\alpha(G) = k - 1$. Applying Theorem~\ref{thm:counting-bound} gives $k - 1 \le \ell(G) \le \alpha(G) = k - 1$, forcing $\ell(G) = k - 1$. If $m\ge1$, Definition~\ref{def:counting-bound} gives
$\ell(G)\le n-1$. If $m=0$, then
$\ell(G)=n=\alpha(G)$.
\end{proof}

\section{Comparison with other lower bounds}

In this section we compare $\ell(G)$ to well-known lower bounds for the independence number, including the Tur\'an-Spencer lower bound, the Caro-Tuza lower bound and the Csaba-Plick-Shokoufandeh bound. In particular, we exhibit infinte families of $k$-uniform hypergraphs for which $\ell(G)$ outperforms these bounds.

\begin{rem}
Let $G$ be a finite simple graph with $n\ge2$ vertices and $m$
edges. Define the Tur\'an lower bound \cite{Turan1941} by
\[
{\rm T}(G)
=
\left\lceil\frac{n^2}{n+2m}\right\rceil.
\]
Equivalently, for $r\in\{2,\ldots,n\}$,
${\rm T}(G)=r$ if and only if
\begin{equation}\label{eq:turan-graph-interval}
\frac{n(n-r)}{2r}
\le m
<
\frac{n(n-r+1)}{2(r-1)}.
\end{equation}
Moreover, ${\rm T}(G)=1$ if and only if
$m=\binom{n}{2}$. For comparison, if $m\ge1$, then for
$r\in\{2,\ldots,n-1\}$,
$\ell(G)=r$ if and only if
\begin{equation}\label{eq:lcnt-graph-interval}
\frac{n(n-1)}{r(r+1)}
\le m
<
\frac{n(n-1)}{r(r-1)}.
\end{equation}
The endpoint cases are
\[
\ell(G)=1
\iff
m=\binom{n}{2}
\]
and
\[
\ell(G)=n
\iff
m=0.
\]
We claim that
\[
\ell(G)\le {\rm T}(G)
\]
for every finite simple graph $G$. Indeed, suppose first that
${\rm T}(G)=r$ for some $2\le r\le n-2$. Then
\[
m\ge\frac{n(n-r)}{2r}.
\]
Furthermore,
\[
\frac{n(n-r)}{2r}
-
\frac{n(n-1)}{r(r+1)}
=
\frac{n(r-1)(n-r-2)}{2r(r+1)}
\ge0.
\]
It follows that
\[
m\ge\frac{n(n-1)}{r(r+1)},
\]
and hence $\ell(G)\le r={\rm T}(G)$. If ${\rm T}(G)=n-1$, then $m\ge1$, while
\[
\frac{n(n-1)}{(n-1)n}=1,
\]
so again $\ell(G)\le n-1$. If ${\rm T}(G)=n$,
then $m=0$ and both bounds equal $n$. Finally, if
${\rm T}(G)=1$, then $G$ is complete and both bounds equal $1$. Therefore, $\ell$ does not improve the classical Tur\'an lower bound in the case of graphs.
\end{rem}

\begin{rem}
For a simple $k$-uniform hypergraph $G$ with $n$ vertices and $m$ edges such that $m\ge \frac{n}{k}$, we define the {\it Tur\'an-Spencer lower bound} to be the number
\begin{equation}
{\rm TS}(G)=\left\lceil \frac{k-1}{k}n\sqrt[k-1]{\frac{n}{km}} \right\rceil\,.
\end{equation}
Assuming $m \ge \frac{n}{k}$, the characterization ${\rm TS}(G)=r$ holds explicitly across all valid ranges of $r$ as follows. For $r \in \{2, \ldots, n\}$, ${\rm TS}(G) = r$ if and only if
    \begin{equation}\label{eq:10}
    \left(\frac{n}{k}\right)^k \left(\frac{k-1}{r}\right)^{k-1} \le m < \left(\frac{n}{k}\right)^k \left(\frac{k-1}{r-1}\right)^{k-1}.
    \end{equation}
For $r = 1$, ${\rm TS}(G) = 1$ if and only if $m \ge \left(\frac{n}{k}\right)^k (k-1)^{k-1}$.
In analogy with the observation made in Remark \ref{rem:5}, the lower bound in \eqref{eq:10} is equivalent to ${\rm TS}(G)\le r$, while the upper bound is equivalent to ${\rm TS}(G)\ge r$. It is well known \cites{Spencer1972, Eustis2013} that ${\rm TS}(G)$ is a valid lower bound for the independence number of $G$ whenever $m \ge \frac{n}{k}$. Direct comparison shows that if $G$ is a graph ($k=2$), then ${\rm TS}(G)\le {\rm T}(G)$ in the range $m\ge \frac{n}{2}$.
\end{rem}

\begin{example}\label{ex:9}
Let $k\ge2$ and suppose that $m$ is an integer in the range
\begin{equation}\label{eq:12}
\frac{n^k}{k^k}\le m < \binom{n}{k}
\end{equation}
Since $n\ge k$, the lower endpoint implies $m\ge n/k$.
Substitution in the definition of ${\rm TS}(G)$ gives
${\rm TS}(G)\le k-1$, while $m<\binom{n}{k}$ and
Proposition~\ref{prop:5} give $\ell(G)\ge k$.
Thus ${\rm TS}(G)<k\le\ell(G)$ for all $k\ge2$ in this range.
Informally, this shows that for any fixed $k \ge 2$ and sufficiently large $n$, $\ell(G)$ strictly outperforms ${\rm TS}(G)$ whenever the hypergraph has a sufficiently high, but not necessarily maximal, edge density. The lower density threshold
$n^k/(k^k\binom{n}{k})$ tends to $k!/k^k$.
In particular, any sequence whose edge density converges to
a value strictly between $k!/k^k$ and $1$ satisfies
\eqref{eq:12} for all sufficiently large $n$.

To see that the range in \eqref{eq:12} contains valid integer edge counts, assume $k \ge 2$. The length of the interval is given by:
\begin{equation*}
\binom{n}{k} - \frac{n^k}{k^k} = n^k \left( \frac{1}{k!} - \frac{1}{k^k} \right) + O(n^{k-1})\,.
\end{equation*}
Because $k^k > k!$ for all $k \ge 2$, the leading coefficient $\left( \frac{1}{k!} - \frac{1}{k^k} \right)$ is strictly positive. Consequently, as $n \to \infty$, the length of the interval grows without bound. Thus, for any fixed $k \ge 2$, the interval has length strictly greater than 1 for all sufficiently large $n$, guaranteeing it contains at least one integer $m$. In particular, for every $k \ge 2$, $\ell(G)$ (and thus, by \eqref{eq:dC}, $\ell_{\rm dC}(G)$) outperforms the Tur\'an-Spencer lower bound on infinitely many $k$-uniform hypergraphs. 
\end{example}

\begin{example}\label{ex:10}
The argument used in Example \ref{ex:9} can be applied in greater generality. For instance, assume $n>6$ and $k=3$. Substituting $r=5$ into \eqref{eq:7}, and substituting $r=4$ into \eqref{eq:10} we see that when the number of edges is in the range
\begin{equation}
\frac{n^3}{108}\le m < \frac{n(n-1)(n-2)}{60}
\end{equation}
then ${\rm TS}(G)<5\le\ell(G)$. Since \[
B_{\rm dC}(n,s,3)\ge B(n,s,3),
\]
for $s\in \{3,4,5\}$ and $n\ge 7$, we also conclude that ${\rm TS}(G)<5\le\ell_{\rm dC}(G)$. 
\end{example}

\begin{rem}\label{rem:11}

For a real number $x$ and a non-negative integer $d$, we use the generalized binomial coefficient defined by the falling product:
\begin{equation*}
\binom{x}{d} = \begin{cases} 
1 & \text{if } d = 0, \\
\frac{x(x-1)\cdots(x-d+1)}{d!} & \text{if } d \ge 1.
\end{cases}
\end{equation*}
Let $G = (V, E)$ be a finite, simple $k$-uniform hypergraph with degree sequence $(d(v))_{v \in V}$. Following Caro and Tuza \cite{CaroTuza1991}, we define the {\it Caro-Tuza lower bound} ${\rm CT}(G)$ for the independence number of $G$ as:
\begin{equation}\label{eq:14}
{\rm CT}(G) = \left\lceil \sum_{v \in V} \binom{d(v)+\frac{1}{k-1}}{d(v)}^{-1} \right\rceil.
\end{equation}
It is shown in \cite{CaroTuza1991} that the Caro-Tuza lower bound is indeed a lower bound for the independence number of finite, simple $k$-uniform hypergraphs. 
\end{rem}

\begin{example} 
Let the vertex set be $V = \mathbb{Z}_6$. First, we construct a $3$-regular $3$-uniform hypergraph $H_1$ using the $6$ cyclic shifts of the edge $\{0,1,3\} \pmod 6$. Every vertex in $H_1$ has degree $3$. The complete $3$-uniform hypergraph $K_6^3$ has $\binom{6}{3} = 20$ edges and is $\binom{5}{2} = 10$-regular. Let $G_1$ be the complement of $H_1$ in $K_6^3$ (i.e., its edge set is $E(K_6^3) \setminus E(H_1)$). Then $G_1$ is a simple $3$-uniform hypergraph with $m = 20 - 6 = 14$ edges, and it is $(10 - 3) = 7$-regular. Direct substitution into \eqref{eq:14} shows that
${\rm CT}(G_1)=2$. On the other hand, substituting $m=14$ into \eqref{eq:7} shows that $\ell(G_1)=3$. 

Similarly, let $H_2$ be the complete bipartite graph $K_{3,3}$ on $6$ vertices. $H_2$ is a $3$-regular graph (i.e., $2$-uniform hypergraph) with $9$ edges. We construct a $4$-uniform hypergraph $G_2$ on the same vertex set by taking the complement of each edge of $H_2$ with respect to the vertex set (i.e., $E(G_2) = \{ V \setminus e \mid e \in E(H_2) \}$). Since every vertex in $H_2$ belongs to exactly $3$ edges, it is excluded from $3$ edges and therefore included in exactly $9 - 3 = 6$ edges of $G_2$. Thus, $G_2$ is a simple $4$-uniform hypergraph with $m = 9$ edges that is $6$-regular. Then ${\rm CT}(G_2)=3$ while $\ell(G_2)=4$. It follows from \eqref{eq:dC} that ${\rm CT}(G_i)< \ell_{\rm dC}(G_i)$ for $i=1,2$.
\end{example}

More generally we have that the separation between $\ell(G)$ and the Caro-Tuza bound can grow linearly with $k$:

\begin{prop}\label{prop:17}
For every integer $k\ge3$, there exist infinitely many
pairwise nonisomorphic finite, simple, regular
$k$-uniform hypergraphs $G$ such that
\[
{\rm CT}(G)\le k-1<k\le\ell(G)\le\ell_{\rm dC}(G)
\]
and
\begin{equation}\label{eq:18bis}
\ell(G) - {\rm CT}(G)>\frac{k-3}{2}\,.
\end{equation}
\end{prop}

\begin{proof}
Fix an integer $k\ge3$. For each positive integer $n$, put
\[
N_n=k(k-1)n,
\qquad
d_n=(kn)^{k-1}.
\]
Since $k$ divides $N_n$, Baranyai's theorem \cite{Baranyai75}
gives a decomposition of the complete $k$-uniform hypergraph
on $N_n$ vertices into $1$-factors. Each $1$-factor is a
collection of $N_n/k$ pairwise disjoint edges covering all
vertices. Hence the number of factors is
\[
F_n
=
\frac{\binom{N_n}{k}}{N_n/k}
=
\binom{N_n-1}{k-1}.
\]
For fixed $k\ge 3$, 
\begin{equation}\label{eq:18}
\frac{F_n}{d_n}
=
\frac{\binom{k(k-1)n-1}{k-1}}{(kn)^{k-1}}
\longrightarrow
\frac{(k-1)^{k-1}}{(k-1)!}
>\frac{(k-1)^{k}}{k!}>1.
\end{equation}
Thus $d_n<F_n$ for all sufficiently large $n$. For each
such $n$, let $G_n$ be the union of any $d_n$ factors in
this decomposition, with the same vertex set. Then $G_n$
is a finite, simple, $d_n$-regular $k$-uniform hypergraph.
Moreover, since at least one factor is omitted, $G_n$ is
not complete. Proposition~\ref{prop:5} therefore gives $\ell(G_n)\ge k$. To evaluate the Caro-Tuza bound, set
\[
t=\frac{1}{k-1}
\]
and denote its unrounded expression by
\[
S_n
=
\frac{N_n}{\binom{d_n+t}{d_n}}.
\]
Since every vertex of $G_n$ has degree $d_n$,
\eqref{eq:14} gives ${\rm CT}(G_n)=\lceil S_n\rceil$.
Using the gamma-function representation of the generalized
binomial coefficient, we obtain
\[
S_n
=
N_n\Gamma(1+t)
\frac{\Gamma(d_n+1)}{\Gamma(d_n+1+t)}.
\]
Stirling's formula implies
\[
\frac{\Gamma(d_n+1)}{\Gamma(d_n+1+t)}
\sim d_n^{-t}.
\]
Since $N_nd_n^{-t}=k-1$, it follows that $S_n\longrightarrow (k-1)\Gamma(1+t)$. Now $0<t<1$, and strict log-convexity of the gamma function
yields
\[
\Gamma(1+t)
<
\Gamma(1)^{1-t}\Gamma(2)^t
=
1.
\]
Consequently, $S_n<k-1$ for all sufficiently large $n$.
Because $k-1$ is an integer, this implies
\begin{equation}\label{eq:19}
{\rm CT}(G_n)
=
\lceil S_n\rceil
\le k-1
<k
\le\ell(G_n).
\end{equation}
Since each selected factor contains $N_n/k$ edges, we have
\[
m_n:=|E(G_n)|=\frac{N_nd_n}{k}\,,
\qquad
\frac{\binom{N_n}{k}}{m_n}=\frac{F_n}{d_n}.
\]
Using the definition of $\ell$ and the arithmetic-geometric mean inequality, we obtain 
\[
\frac{F_n}{d_n}\le \binom{\ell(G_n)+1}{k}\le \frac{\left(\ell(G_n)-\frac{k-3}{2}\right)^k}{k!}\,.
\]
from which, using \eqref{eq:18}, we obtain
\begin{equation}\label{eq:20}
\ell(G_n)>\frac{k-3}{2}+k-1
\end{equation}
for fixed $k$ and sufficiently large $n$. Combining \eqref{eq:19} and \eqref{eq:20}, yields \eqref{eq:18bis}.

It remains to compare $\ell(G_n)$ and $\ell_{\rm dC}(G_n)$.
The number of edges of $G_n$ is
\[
m_n=\frac{N_nd_n}{k}\longrightarrow\infty.
\]
For fixed $k$, we have
\[
B_{\rm dC}(N_n,N_n-k+1,k)
=
\frac{\binom{N_n}{k}}{\binom{N_n-k+1}{k}}
\longrightarrow 1\,.
\]
Thus, for all sufficiently large $n$,
\[
B_{\rm dC}(N_n,N_n-k+1,k)\le m_n\, ,
\]
and the definition gives
\[
\ell_{\rm dC}(G_n)\le N_n-k<N_n-k+1\, .
\]
Applying the conditional comparison in
Remark~\ref{rem:covering-comparison}, we obtain
\[
\ell(G_n)\le\ell_{\rm dC}(G_n)\, .
\]
Together with the preceding Caro-Tuza estimate, this proves
\[
{\rm CT}(G_n)<\ell(G_n)\le\ell_{\rm dC}(G_n)
\]
for all sufficiently large $n$.
Finally, the orders $N_n=k(k-1)n$ are strictly increasing,
so these hypergraphs provide infinitely many pairwise
nonisomorphic examples. 
\end{proof}

\begin{rem}
Let $G$ be a $3$-uniform hypergraph. Borrowing the notation of Remark \ref{rem:11} we define the {\it Csaba-Plick-Shokoufandeh lower bound} to be
\begin{equation}
{\rm CPS}(G)=\left\lceil\frac{\sqrt{\pi}}{2}\sum_{v\in V(G)} \frac{1}{\sqrt{d(v)+1}}\right\rceil\,.
\end{equation}
It is shown in \cite{CsabaPlickShokoufandeh2012} that ${\rm CPS}(G)$ is indeed a lower bound for the independence number of $G$. 
\end{rem}

\begin{example}
Let $G_n$ be obtained from the $3$-uniform complete hypergraph on $n$ vertices by removing a single edge so that ${\rm CPS}(G_n)=\left \lceil g(n) \right \rceil$, where
\begin{equation}
g(n)=\sqrt{\frac{\pi}{2}}\left(\frac{n-3}{\sqrt{n^2-3n+4}}+\frac{3}{\sqrt{n^2-3n+2}} \right)
\end{equation}
Since for real $x\ge3$, the derivative is
\[
g'(x)=\sqrt{\frac\pi2}\left(
\frac{3x-1}{2(x^2-3x+4)^{3/2}}
-\frac{3(2x-3)}{2(x^2-3x+2)^{3/2}}
\right)<0\,,
\]
$g(n)$ is decreasing for $n>3$ and
\[
\lim_{n\to\infty}g(n)=\sqrt{\frac{\pi}{2}}>1\,,
\qquad
g(4)=\sqrt{\pi}\left(\frac14+\frac{\sqrt3}{2}\right)<2\,.
\]
Hence, for every integer $n\ge4$, $1<g(n)\le g(4)<2$
and therefore ${\rm CPS}(G_n)=\lceil g(n)\rceil=2$. On the other hand, Proposition~\ref{prop:5} gives $\ell(G_n)\ge3$.
Thus this family provides infinitely many $3$-uniform hypergraphs
for which $\ell(G)$ (and thus, by \eqref{eq:dC}, $\ell_{dC}(G)$) exceeds the Csaba-Plick-Shokoufandeh bound. It is worth noting that in this example $\ell$ outperforms the Csaba-Plick-Shokoufandeh bound because the latter fails to certify even the obvious bound $\alpha(G)\ge 3$, which holds for all non-complete $3$-uniform hypergraphs. 
\end{example}

\begin{example}
Proposition~\ref{prop:17}. In this case,
\[
N_n=6n,\qquad d_n=9n^2,\qquad
F_n=\binom{6n-1}{2}\, .
\]
Since
\[
F_n-d_n=9n(n-1)+1>0\, ,
\]
the construction gives a noncomplete hypergraph $G_n$
for every positive integer $n$. Moreover,
\[
1<
\frac{3n\sqrt{\pi}}{\sqrt{9n^2+1}}
<2\,.
\]
It follows that, for every positive integer $n$,
\[
{\rm CPS}(G_n)
=
\left\lceil
\frac{3n\sqrt{\pi}}{\sqrt{9n^2+1}}
\right\rceil
=2<3\le\ell(G_n)\,.
\]
\end{example}

\section{Conclusions and further directions}

We have studied the lower bound $\ell(G)$ obtained by inverting
the classical counting inequality for Tur\'an systems. This bound
depends only on the number of vertices, the number of edges, and
the uniformity of the hypergraph. For fixed uniformity, it can be
evaluated exactly by binary search in time polynomial in the binary
input length. The contribution of the present paper concerns the
explicit evaluation and comparison of this bound, rather than a
new covering inequality.

Our comparisons distinguish the graph case from the case of
higher uniformity. For graphs, $\ell(G)$ never exceeds the rounded
average-degree bound following from Tur\'an's theorem. For every
fixed $k\ge3$, however, we identify parameter ranges containing
infinitely many $k$-uniform hypergraphs on which $\ell(G)$ exceeds
the Tur\'an-Spencer bound, and we construct infinite families
on which it exceeds the Caro-Tuza bound. We also obtain a
separation from the Csaba-Plick-Shokoufandeh bound in the
$3$-uniform case. 

Several questions remain concerning the strengths and limitations
of this approach. It would be useful to characterize more fully
the hypergraphs for which $\ell(G)=\alpha(G)$ and to better understand
the parameter regimes in which the counting bound provides
a useful complement to other estimates. Another direction is
to carefully quantify the improvements obtained by inverting stronger
classical bounds for covering and Tur\'an numbers, while keeping
track of the computational cost of evaluating the resulting
independence bounds.

Finally, our original algebraic motivation suggests a broader
direction of investigation. The cohomological approach to
independence developed in \cite{AGGHK}, together with the
association of $2$-step nilpotent $L_\infty$-algebras to
hypergraphs in \cite{AldiBevins}, provides a setting in which
to reconsider classical independence inequalities. It would
be interesting to better understand which lower bounds for
the independence numbers of graphs and hypergraphs admit
generalizations to $2$-step nilpotent $L_\infty$-algebras,
including algebras that do not arise from these combinatorial
constructions. In particular, one may ask which counting,
probabilistic, and degree-based arguments have formulations
in terms of intrinsic algebraic data, and what additional
hypotheses or modifications are needed for such formulations
to remain valid. Identifying both extensions and obstructions
would help clarify which features of classical independence
theory are genuinely combinatorial and which persist in this
broader algebraic setting.

\subsection*{Acknowledgements} The work on this paper was supported by the National Science Foundation grant number DMS1950015. We are very grateful to Craig Larson for illuminating correspondence and encouragement. We would like to thank the anonymous referee for extremely valuable feedback. During the preparation of this work the authors used ChatGPT 6, Claude Opus 5, and Gemini 3.1 in order to improve the exposition. After using this tool, the authors reviewed and edited the content as needed and take full responsibility for the content of the published article.

\end{document}